\documentclass{article}
\usepackage[utf8]{inputenc}
\usepackage{ gensymb }
\usepackage{amsmath}
\usepackage{ amssymb }
\usepackage{amsthm}
\newtheorem{Main}{Theorem}
\newtheorem{Lemma 1}[Main]{Lemma}
\newtheorem{Lemma 2}[Main]{Lemma}
\usepackage{enumitem}
\usepackage{hyperref}
\usepackage{graphicx}
\title{ON A GENERALIZATION OF TUPPER'S FORMULA FOR $m$ COLOURS  AND $n$ DIMENSIONS}
\author{Sai Teja Somu \and 
	Vidyanshu Mishra}
\date{\today}

\begin{document}
	\maketitle
	\begin{abstract}
		Tupper's formula $\frac{1}{2}<\bigg\lfloor \bmod \bigg(\lfloor \frac{y}{17}\rfloor 2^{-17\lfloor x \rfloor -\bmod (\lfloor y \rfloor,17)},2\bigg)\bigg\rfloor$ has an interesting property that for any monochrome image that can be represented by pixels in a two dimensional array of dimensions $106\times 17$, there exists a natural number $k$ such that the graph of the equation in the range $0\leq x <106$ and $k\leq y<k+17$, is that image. In this paper, we give a generalization for $m$ colours and $n$ dimensions. We give $m$ formulae consisting of $n$ free variables, with the property that, for any $n$ dimensional object of $m$ colours $C_1,\cdots, C_m$, that can be represented by hypervoxels(multidimensional analogue of pixel) in a $n$ dimensional array of dimensions $A_1\times \cdots \times A_n$, there exists a natural number $k$ such that, when the first formula is graphed using colour $C_1$, second formula is graphed using colour $C_2$,$\cdots$, $m$th formula is graphed using colour $C_m$ in the range $0\leq x_1<A_1$,$0\leq x_2<A_2,\cdots, 0\leq x_{n-1}<A_{n-1},k\leq x_n <k+A_n$, the union of all graphs is that $n$-dimensional object.
	\end{abstract}
	\section{Introduction}
	In \cite{Tupper}, Tupper introduced a formula \begin{equation}\frac{1}{2}<\bigg\lfloor \bmod \bigg(\big\lfloor \frac{y}{17}\big\rfloor 2^{-17\lfloor x \rfloor-\bmod(\lfloor y \rfloor,17)},2\bigg)\bigg\rfloor,\end{equation} which has an interesting property that for any monochrome image that can be represented in pixels in a two dimensional array of dimensions $106\times 17$, there exists a natural number $k$ such that the graph of the equation for $0\leq x<106$ and $k\leq y <k+17$ is that image( see Lemma 1 of \cite{Trevino} for a proof). Since the formula given by Tupper can be drawn in a $106\times 17$ grid, there exists a value of $k$  for which the graph of the formula in $0\leq x<106$ and $k\leq y <k+17$ is an image of Tupper's formula. While graphing (1) in Figure 1, we are painting all ordered pairs $(x,y)$ that satisfy the inequality with black. 
	
	\begin{figure}[htp]
		\centering
		\includegraphics[width=6cm]{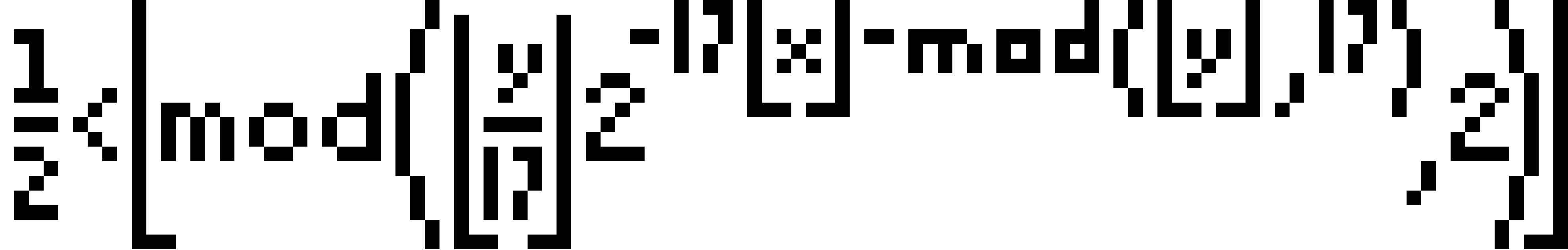}
		\caption{Graph of (1) in the range $0\leq x<106$ and $k\leq y <k+17$ for a known value of $k$}
		\label{fig:tupper}
	\end{figure}
	Now a question arises, can we generalize this property for images consisting of multiple colours by considering multiple formulae? There is a generalization. Suppose we have to graph multiple formulae (2), (3), (4),
	\begin{align}
	&\frac{1}{2}<\bigg\lfloor \bmod\bigg(\left\lfloor\frac{y}{68} \right\rfloor2^{-68^2\lfloor x \rfloor-68\bmod(\lfloor y\rfloor,68)-1},2\bigg) \bigg\rfloor,\\
	&\frac{1}{2}<\bigg\lfloor \bmod\bigg(\left\lfloor\frac{y}{68} \right\rfloor2^{-68^2\lfloor x \rfloor-68\bmod(\lfloor y\rfloor,68)-2},2\bigg) \bigg\rfloor,\\
	&\frac{1}{2}<\bigg\lfloor \bmod\bigg(\left\lfloor\frac{y}{68} \right\rfloor2^{-68^2\lfloor x \rfloor-68\bmod(\lfloor y\rfloor,68)-3},2\bigg) \bigg\rfloor,
	\end{align}
	in the range $0\leq x<50$ and $k\leq y<k+15$ such that $(x,y)$ satisfying inequality (2) will be painted with the colour blue, $(x,y)$ satisfying (3) will be painted with the colour red, and $(x,y)$ satisfying inequality (4) will be painted with the colour green then for any image consisting of three colours red, blue, green, there exists a natural number $k$ such that the union of the three graphs (2), (3) and (4)  will be that image. For instance, there exist natural numbers $k_1,k_2$ such that union of graphs of (2), (3), (4) in $0\leq x<50$ and $k_1\leq y<k_1+15$ is Figure 2 and union of graphs of (2), (3), (4) in $0\leq x<50$ and $k_2\leq y<k_2+15$ is Figure 3.   
	
	\begin{figure}[htp]
		\centering
		\includegraphics[width=6cm]{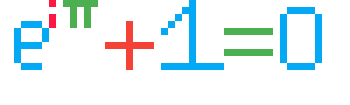}
		\caption{Graph of (2) in blue, (3) in red, (4) in green in the range $0\leq x<50$ and $k_1\leq y <k_1+15$ for a determinable value of $k_1$}
		\label{fig:euler}
	\end{figure}
	
	\begin{figure}[htp]
		\centering
		\includegraphics[width=6cm]{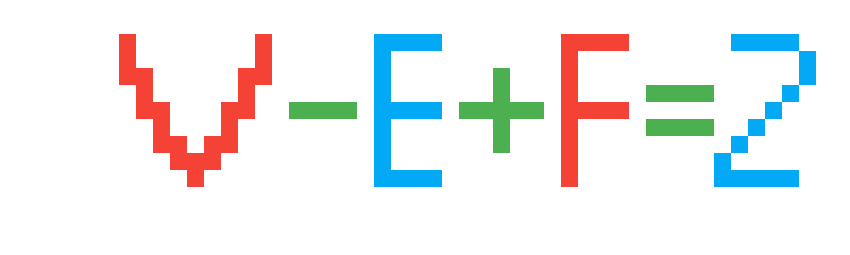}
		\caption{Graph of (2) in blue, (3) in red, (4) in green in the range $0\leq x<50$ and $k_2\leq y <k_2+15$ for a determinable value of $k_2$}
		\label{fig:polyhedron}
	\end{figure}
	
	We prove this generalization for $m$ colours, and $n$ dimensions in the form of this theorem.
	
	\begin{Main}
		For any given positive integers $n\geq 2,m,A_1,A_2,\cdots,A_n$ let $R=A_1+\cdots +A_n+m$ and for $1\leq i \leq m$ define $$f_i(x_1,\cdots,x_n)=\bigg\lfloor\bmod\bigg(\bigg\lfloor\frac{x_n}{R}\bigg\rfloor2^{-R^n\lfloor x_1\rfloor-R^{n-1}\lfloor x_2\rfloor-\cdots-R^2\lfloor x_{n-1}\rfloor-R\bmod(\lfloor x_n\rfloor,R)-i},2\bigg) \bigg\rfloor,$$
		then for any $m$ pairwise disjoint subsets $S_1,\cdots, S_m$ of $S=[0,A_1)\times\cdots\times [0,A_n)\cap \mathbb{Z}^n,$ there exists a natural number $k$ given by 
		\begin{equation}
		k= R\sum_{j=1}^{m}\sum_{(m_1,\cdots,m_n)\in S_j}2^{R^nm_1+R^{n-1}m_{2}+\cdots + Rm_n+j}
		\end{equation} 
		such that \begin{multline*}
		\{(x_1,x_2,\cdots,x_n):\frac{1}{2}<f_i(x_1,\cdots,x_n)\text{ and }(x_1,\cdots,x_{n-1},x_n-k)\in[0,A_1)\times\cdots\times [0,A_n)\}\\= \{(x_1,x_2,\cdots,x_n):(\lfloor x_1\rfloor,\lfloor x_2\rfloor,\cdots,\lfloor x_{n-1}\rfloor,\lfloor x_n-k\rfloor)\in S_i\},
		\end{multline*}   for all
		$1\leq i \leq m$.\end{Main}
	Formulae (2), (3), and (4) are special cases of Theorem 1 when $A_1=50,$ $A_2=15$, $m=3$, $n=2$. For any $n$ dimensional object of $m$ colours $C_1,\cdots, C_m$, that can be represented by hypervoxels in a $n$ dimensional array of dimensions $A_1\times \cdots \times A_n$, let $S_1$ be the set of tuples corresponding to colour $C_1$ for that object, $\cdots$, let $S_m$ be the set of tuples corresponding to colour $C_m$ for that object. Then for the value of $k$ given by (5),   when the graph of $\frac{1}{2}<f_1$ is coloured using colour $C_1$,  the graph of $\frac{1}{2}<f_2$ is coloured using colour $C_2$, $\cdots$, the graph of $\frac{1}{2}<f_m$ is coloured using colour $C_m$ in the range $0\leq x_1<A_1,\cdots,0\leq x_{n-1}<A_{n-1},k\leq x_n<k+A_n$ we get that $n$-dimensional object. 
	\section{Proof of Theorem 1}
	We will prove Theorem 1 using some lemmas.
	\begin{Lemma 1}
		For $k$ given by (5), if $x_n-k\in [0,A_n)$ then 
		$$\bigg\lfloor\frac{x_n}{R} \bigg\rfloor=\sum_{j=1}^{m}\sum_{(m_1,\cdots,m_n)\in S_j}2^{R^nm_1+R^{n-1}m_{2}+\cdots + Rm_n+j},$$ and $$\bmod(\lfloor x_n\rfloor,R) = \lfloor x_n-k \rfloor.$$
	\end{Lemma 1}
	\begin{proof}
		As $R=A_1+\cdots+A_n+m$, we have $R>A_n$ and as $x_n-k\in [0,A_n)$ we have $0\leq \frac{x_n-k}{R}<1$. As $\frac{x_n}{R}=\frac{k}{R}+\frac{x_n-k}{R}$, as $\frac{k}{R}\in \mathbb{Z}$ and $\frac{x_n-k}{R}\in [0,1)$ we have $$\bigg\lfloor\frac{x_n}{R} \bigg\rfloor=\frac{k}{R}= \sum_{j=1}^{m}\sum_{(m_1,\cdots,m_n)\in S_j}2^{R^nm_1+R^{n-1}m_{2}+\cdots + Rm_n+j}.$$
		For the second part as $\lfloor x_n \rfloor = \lfloor x_n-k \rfloor+k$, and $k$ is a multiple of $R$ and $0\leq \lfloor x_n-k \rfloor<A_n<R$, we have $$\bmod(\lfloor x_n \rfloor,R)=\lfloor x_n-k\rfloor.$$
	\end{proof}
	\begin{Lemma 2}
		If $\alpha = \sum_{i=1}^{k}2^{n_i}$ for distinct integers $n_1,\cdots,n_k$ then $$\frac{1}{2}<\lfloor\bmod(\alpha ,2)\rfloor\text{ if and only if } 0\in \{n_1,\cdots,n_k\}.$$
	\end{Lemma 2}
	\begin{proof}
		Let $N$ be the subset of $A=\{n_1,\cdots,n_k\}$ consisting of all negative integers of $A$, $P$ be the subset of $A$ consisting of all positive integers of $A$. Note that $$\alpha = \sum_{a\in N}2^{a}+\sum_{a\in P}2^{a}+\epsilon,$$ where $\epsilon=1$, if $0\in A$ and $\epsilon=0$ if $0\notin A$. As $\sum_{a\in P}2^{a}$ is even and $0\leq \sum_{a\in N}2^{a}<1$  we have $\bmod(\alpha,2)=\sum_{a\in N}2^{a}+\epsilon$. Now as, $0\leq \sum_{a\in N}2^{a}<1$ we have $\lfloor\bmod(\alpha,2)\rfloor=\epsilon$. 
		Therefore, $\lfloor\bmod(\alpha,2)\rfloor=1>\frac{1}{2}$ if $0\in \{n_1,\cdots,n_k\}$ and $\lfloor\bmod(\alpha,2)\rfloor=0\leq \frac{1}{2}$ if $0\notin \{n_1,\cdots,n_k\}$ which implies $\frac{1}{2}<\lfloor\bmod(\alpha,2)\rfloor$ if and only if $0\in \{n_1,\cdots,n_k\}$. 
	\end{proof}
	Now we are ready to prove Theorem 1.
	\begin{proof}
		Let $(x_1,\cdots,x_{n-1},x_{n})$ be a tuple such that $(x_1,\cdots,x_{n-1},x_{n}-k)\in [0,A_1)\times \cdots \times [0,A_n)$. As $x_n-k\in [0,A_n)$, from Lemma 2, we have \begin{equation}\bigg\lfloor\frac{x_n}{R} \bigg\rfloor=\sum_{j=1}^{m}\sum_{(m_1,\cdots,m_n)\in S_j}2^{R^nm_1+R^{n-1}m_{2}+\cdots + Rm_n+j},\end{equation} and \begin{equation}\bmod(\lfloor x_n\rfloor,R) = \lfloor x_n-k \rfloor.\end{equation} Now using (6) and (7) for any $1\leq i \leq m$ we have  \begin{align*}
		f_i(x_1,\cdots,x_n)&=\bigg\lfloor\bmod\bigg(\bigg\lfloor\frac{x_n}{R}\bigg\rfloor2^{-R^n\lfloor x_1\rfloor-R^{n-1}\lfloor x_2\rfloor-\cdots-R^2\lfloor x_{n-1}\rfloor-R\bmod(\lfloor x_n\rfloor,R)-i},2\bigg) \bigg\rfloor,\\
		&=\bigg\lfloor\bmod\bigg(\bigg\lfloor\frac{x_n}{R}\bigg\rfloor2^{-R^n\lfloor x_1\rfloor-R^{n-1}\lfloor x_2\rfloor-\cdots-R^2\lfloor x_{n-1}\rfloor-R\lfloor x_n-k \rfloor-i},2\bigg)\bigg\rfloor,\\
		&=\bigg\lfloor\bmod\bigg(\sum_{j=1}^{m}\sum_{(m_1,\cdots,m_n)\in S_j}2^{R^n(m_1-\lfloor x_1\rfloor)+R^{n-1}(m_{2}-\lfloor x_2\rfloor)+\cdots + R(m_n-\lfloor x_n-k\rfloor)+(j-i)},2\bigg)\bigg\rfloor.
		\end{align*}
		Therefore, $$f_i(x_1,\cdots,x_n)=\bigg\lfloor\bmod\bigg(\sum_{j=1}^{m}\sum_{(m_1,\cdots,m_n)\in S_j}2^{R^n(m_1-\lfloor x_1\rfloor)+R^{n-1}(m_{2}-\lfloor x_2\rfloor)+\cdots + R(m_n-\lfloor x_n-k\rfloor)+(j-i)},2\bigg)\bigg\rfloor.$$ Notice that in the above equation $R^n(m_1-\lfloor x_1\rfloor)+R^{n-1}(m_{2}-\lfloor x_2\rfloor)+\cdots + R(m_n-\lfloor x_n-k\rfloor)+(j-i)$ are distinct integers for different tuples $(m_1,m_2,\cdots, m_n)\in S_j$ as $R>m_1,m_2,\cdots, m_n,j$. So from Lemma 3, $\frac{1}{2}<f_i(x_1,\cdots,x_n)$ if and only if $R^n(m_1-\lfloor x_1\rfloor)+R^{n-1}(m_{2}-\lfloor x_2\rfloor)+\cdots + R(m_n-\lfloor x_n-k\rfloor)+(j-i)=0$ for some $(m_1,\cdots,m_n)\in S_j$ and as $R^n(m_1-\lfloor x_1\rfloor)+R^{n-1}(m_{2}-\lfloor x_2\rfloor)+\cdots + R(m_n-\lfloor x_n-k\rfloor)+(j-i)=0$ if and only if $\lfloor x_1\rfloor=m_1$, $\lfloor x_2\rfloor=m_2$,$\cdots$,$\lfloor x_{n-1}\rfloor=m_{n-1}$ and $\lfloor x_n-k \rfloor=m_n$, $j=i$ we have $\frac{1}{2}<f_i(x_1,\cdots,x_n)$ if and only if $(\lfloor x_1\rfloor,\cdots,\lfloor x_{n-1} \rfloor,\lfloor x_n-k \rfloor)\in S_i$. Therefore, for all $1\leq i \leq m$ we have, \begin{multline*}
		\{(x_1,x_2,\cdots,x_n):\frac{1}{2}<f_i(x_1,\cdots,x_n)\text{ and }(x_1,\cdots,x_{n-1},x_n-k)\in[0,A_1)\times\cdots\times [0,A_n)\}\\= \{(x_1,x_2,\cdots,x_n):(\lfloor x_1\rfloor,\lfloor x_2\rfloor,\cdots,\lfloor x_{n-1}\rfloor,\lfloor x_n-k\rfloor)\in S_i\}.
		\end{multline*}   
	\end{proof}
	\bibliographystyle{amsplain}

\end{document}